\documentclass{article}
\usepackage[english]{babel}
\usepackage{amssymb}
\usepackage{amsthm}
\usepackage{pstricks}
\usepackage{epic,eepic,graphicx}
\usepackage{amsmath}

\newtheorem{theorem}{Theorem}[section]
\newtheorem{lemma}[theorem]{Lemma}
\newtheorem{proposition}[theorem]{Proposition}
\newtheorem{corollary}[theorem]{Corollary}
\newtheorem{conjecture}[theorem]{Conjecture}

{\theoremstyle{definition}
\newtheorem{definition}[theorem]{Definition}
\newtheorem{remark}[theorem]{Remark}}

\newcommand{\bt}{\begin{theorem}}
\newcommand{\et}{\end{theorem}}
\newcommand{\bl}{\begin{lemma}}
\newcommand{\el}{\end{lemma}}
\newcommand{\bpf}{\begin{proof}}
\newcommand{\epf}{\end{proof}}
\newcommand{\bp}{\begin{proposition}}
\newcommand{\ep}{\end{proposition}}
\newcommand{\bc}{\begin{corollary}}
\newcommand{\ec}{\end{corollary}}
\newcommand{\bcc}{\begin{conjecture}}
\newcommand{\ecc}{\end{conjecture}}
\newcommand{\bd}{\begin{definition}}
\newcommand{\ed}{\end{definition}}
\newcommand{\br}{\begin{remark}}
\newcommand{\er}{\end{remark}}

\def\bysame{\leavevmode\hbox to3em{\hrulefill}\thinspace}

\begin{document}

\title{{\huge 50 years of Finite Geometry, \\the ``geometries over finite rings" part}}

\author{{\Large Dirk Keppens} \bigskip \\Faculty of Engineering Technology, KU Leuven\\Gebr. Desmetstraat 1 \\B-9000 Ghent 
BELGIUM \\\textbf{\texttt{dirk.keppens@kuleuven.be}}}




\date{}


\maketitle

\begin{abstract}
Whereas for a substantial part, ``Finite Geometry"  during the past 50 years has focussed on geometries over finite fields, geometries over finite rings that are not division rings have got less attention. Nevertheless, \mbox{several} important classes of finite rings give rise to interesting geometries.\newline In this paper we bring together some results, scattered over the literature, con\-cerning finite rings and plane projective geometry over such rings. It doesn't contain new material, but by collecting stuff in one place, we hope to stimulate further research in this area for at least another 50 years of Finite Geometry. 
\end{abstract}

\bigskip


\noindent


\section{Introduction}

Geometries over rings that are not division rings have been studied for a long time.
The first systematic study was done by Dan Barbilian \cite{Barb}, besides a mathe\-matician also one of the greatest Romanian poets (with pseudonym Ion Barbu). He introduced plane projective geometries over a class of associative rings with unit, called Z--rings (abbreviation for {\it Zweiseitig singul\"are Ringe}) which today are also known as Dedekind-finite rings. These are rings with the property that $ab=1$ implies $ba=1$ and they include of course all commutative rings but also all finite rings (even non--commutative). 

Wilhelm Klingenberg introduced in \cite{Kling4} projective planes and 3--spaces over {\it local rings}. A ring $R$ is local if it possesses a unique maximal right ideal (which turns out to be the Jacobson radical $J(R)$). 
For a local ring $R$ the quotient ring $R/J(R)$ is a division ring (= skewfield or field) and the natural homomorphism of $R$  onto $\Bbb{K}=R/J(R)$ induces an epimorphism of the plane $P_2(R)$ over $R$ onto the ordinary desarguesian projective plane PG(2,$\Bbb{K}$). Nowadays planes over local rings are called (desarguesian) {\it Klingenberg planes} (see also \cite{Bacon}). In the finite case such planes have the finite projective plane PG(2,$q$) over the Galois field GF($q$) as epimorphic image. 

In three other papers \cite{Kling1}, \cite{Kling2} and \cite{Kling3}, Klingenberg studied projective planes over local rings with some additional properties, called $H$--rings (short for {\it Hjelmslev rings}). In these rings  the left and right ideals form a chain and the maximal ideal contains only zero divisors. If one drops that last condition, one gets {\it chain rings}. In the finite case any chain ring is an $H$--ring.
Planes over $H$--rings are now called (desarguesian) {\it Hjelmslev planes} after the Danish mathe\-matician Johannes Hjelmslev (born as Johannes Petersen) who was the first one to consider plane geometries in which two distinct lines may have more than one point in common \cite{Hj}. Among the finite $H$--rings are the Galois rings GR($p^{nr},p^n$) of cardinality $p^{nr}$ and characteristic $p^n$ which are natural generali\-zations of Galois fields. 

In the early seventies another class of rings came under the attention: {\it full matrix rings} over fields. Strongly inspired by the work of the italian ``father of Galois geometry" Beniamino Segre on geometries over finite fields (e.g.~\cite{BSegre}), J.A. Thas defined projective planes (and higher dimensional spaces) over full matrix rings with elements in a field and investigated combinatorial properties in the finite planes over the matrix rings $M_n(GF(q))$ of $n \times n$--matrices over Galois fields \cite{Thas}. We will refer to these planes further as {\it Thas planes}.

In the eighties F.D. Veldkamp was very productive in the area of projective ring planes and their generalizations. He gave in \cite{Veld1} and \cite{Veld3} an axiomatic description of projective planes and higher dimensional geometries over the large class of {\it rings of stable rank 2}, a notion coming from algebraic $K$--theory. \newline A ring $R$ has stable rank 2 if for any $a,b \in R$ with $Ra+Rb=R$ there exists $r \in R$ such that $a+rb$ is a unit. The class of rings of stable rank 2 includes the class of semilocal rings (hence also all finite rings, local rings, chain rings, $H$--rings and matrix rings over a division ring) and a ring of stable rank 2 is always Dedekind--finite (hence a Z--ring in the sense of Barbilian). Projective planes over rings of stable rank 2 are called (desarguesian) {\it Veldkamp planes}. Among these are Klingenberg planes, Hjelmslev planes, Thas planes and also the projective planes over semiprimary rings (i.e.~rings with nilpotent Jacobson radical and with $R/J(R)$ semisimple) treated by Bingen in \cite{Bingen}.

In almost all papers on projective geometry over rings no special attention is paid to the finite case. Mostly, theorems deal with rings in general (with no specification for finite or infinite).
In this paper we restrict ourselves to the finite case. First we bring together some results on finite rings with special attention for local rings. Then we have a closer look at projective plane geometries over finite rings. In the last section we deal with some generalizations of rings (semi\-rings, nearrings and alternative rings) and projective plane geometries over such algebraic structures.

\section{Finite rings}

In this section the word ``ring"  always refers to an associative ring with unit $1 \not= 0$, but with multiplication not necessarily commutative.

Finite fields or Galois fields are well-known algebraic structures. Finite fields of order $q$ only exist if $q$ is a prime power ($q=p^r$) and for each such $q$ their is a unique (up to isomorphism) field of that order which is denoted by $\Bbb{F}_q$ or by GF($q$). The prime number $p$ is the characteristic of the field. 

It is natural to look at generalizations of finite fields to finite rings, but the situation is much more complicated. First there exist finite non--commutative rings unlike the situation in finite division rings where the famous theorem of Wedderburn forces any finite skewfield to be a field. Also the order of a finite ring doesn't uniquely determine that ring (for example there are four non--isomorphic rings of order four, including one field).
A complete classification of finite rings seems to be a ``mission impossible" (even if one restricts to the commutative case).\newline
The paper of Raghavendran \cite{Rag} on rings of prime power order was the starting point for the study of the structure of finite rings. Also the work of Wilson \cite{Wilson1} and \cite{Wilson2} was of great importance. A recent survey on results obtained so far with an extensive bibliography can be found in Nechaev \cite{Nech}. 

Local rings, first defined by Krull in \cite{Krull}, play a central role in the structure theory of (finite) rings. 
Recall that a ring $R$ is called {\it local} if it possesses a unique maximal right ideal (or equivalently a unique maximal left ideal). This is stronger than asking that $R$ has a unique maximal two--sided ideal (e.g.~the ring $M_n(\Bbb{Z}_{p^n})$ of $n \times n$--matrices over $\Bbb{Z}/p^n\Bbb{Z}$ has a unique maximal two--sided ideal but is not local). The unique maximal right or left ideal in a local ring turns out to be the Jacobson radical $J(R)$. Other characterizations of local rings are possible. E.g.~$R$ is local iff the set of non--units forms a right (or left) proper ideal in $R$. Also $R$ is local if and only if for all $r \in R$ either $r$ or $1-r$ is invertible. Finally $R$ is local iff $R/J(R)$ is a division ring. Other characterizations in terms of zero divisors are given in \cite{Sat}.
In the finite case one can say even more: $R$ is local iff $R \setminus J(R)$ is the set of units of $R$ or equivalently $J(R)$ is the set of nilpotent elements of $R$. 

The following theorem gives parameters for finite local rings.

\begin{theorem}{\rm (Raghavendran \cite{Rag})}
Let $R$ be a finite local ring. Then there exist unique numbers $p$, $n$, $r$ and $k$ such that $|R|=p^{nr}$, $|J(R)|=p^{(n-1)r}$ and the characteristic of $R$ is $p^k$ with $1 \leq k \leq n$. The number $p^r$ is the order of the Galois field $R/J(R)$ and the number $n$ is the index of nilpotency of $J(R)$. If $k=n$ then $R$ is commutative.

\end{theorem}

There is also a more recent result which conversely characterizes local rings among finite rings just by a couple of parameters.

\begin{theorem}{\rm (Behboodi and Beyranvand \cite{BeBe}, Gonz\'alez \cite{Gonz})}
Let $R$ be a finite ring and let $Z(R)$ be the set of zero--divisors of $R$. Then $R$ is local if and only if $|R|=p^m$ and $|Z(R)|=p^n$ for some prime number $p$ and integers $1 \leq n < m$. Moreover, when $R$ is local with these parameters, then the order of $R/J(R)=R/Z(R)=p^r$ with $r=m-n$.
\end{theorem}

The structure of commutative finite local rings was first studied by Ganske and McDonald \cite{GaMcD}. Classification theorems are proved for fixed orders or fixed characteristic in \cite{Chi}, \cite{Chi2}, \cite{Corb} and \cite{Rag}. In the non--commutative case Wirt \cite{Wirt} has contributed to the theory.

By the following important structure theorem the classification problem of commutative finite rings can be reduced to that of finite local rings.
\begin{theorem}{\rm (McDonald \cite{McD})}
	Let $R$ be a finite commutative ring. Then $R$ decomposes (up to order of summands) uniquely as a direct sum of finite local rings.
\end{theorem}

Another decomposition theorem, also valid in the non--commutative case, shows once more the importance of finite local rings.

\begin{theorem}{\rm (McDonald \cite{McD} and Wirt \cite{Wirt})}
	Let $R$ be a finite ring. Then $R$ decomposes as $S+N$ with $S$ a direct sum of full matrix rings over finite local rings and $N$ a subring of the Jacobson radical $J(R)$.
\end{theorem}

Next we look at {\it principal ideal rings}. A ring is called a right principal ideal ring if any right ideal $I$ is right principal, i.e. generated by one element ($I=aR$). Similar definition for left principal ideal ring. If a ring is a left and right principal ideal ring, it is called a principal ideal ring (PIR).
A right principal ideal ring is always a right noetherian ring, since any right ideal is finitely generated. It is also a right B\'ezout ring since any finitely generated right ideal is principal. In fact the  right PIR's are just the rings which are both right noetherian and right B\'ezout. Similar results are true for left principal ideal rings and PIR's.\newline
The structure of finite principal ideal rings was first studied by Fisher in \cite{Fisher}. 
For finite rings the notions of left PIR, right PIR and PIR are equivalent (see \cite{Nech}).

Another important class of rings are the {\it chain rings}. A ring is called a right chain ring if for any $a$ and $b$ in $R$ either $a \in bR$ or $b \in aR$. For a right chain ring the lattice of right ideals is totally ordered by inclusion and it follows that $R$ is a local ring. Analogous definitions and results for left chain rings. A ring which is a left and right chain ring is called a chain ring. In the infinite case there are examples of right chain rings which are not chain rings (see \cite{LoLa}, \cite{Sk} and \cite{BrTo}), but in the finite case there is a left--right equivalence.
Every ideal of a chain ring is a power of the unique maximal ideal. A (left or right) chain ring with the additional property that any non--unit is a two--sided zero divisor is called a (left or right) $H$--ring or {\it affine Hjelmslev ring} (two--sided chain rings are also known as {\it projective Hjelmslev rings}). Finite chain rings are always left and right $H$--rings. For a comprehensive study of $H$--rings linked to the behaviour of ideals, we refer to \cite{To}.

The following theorem shows that finite chain rings are nothing but finite local principal ideal rings!

\begin{theorem}{\rm (Clarke and Drake \cite{ClD} and Lorimer \cite{Lorim3})}
Let $R$ be a finite ring. Then the following conditions are equivalent:
(a) $R$ is a local PIR
(b) $R$ is a local ring with principal maximal ideal
(b) $R$ is a left chain ring
(c) $R$ is a right chain ring
(d) $R$ is a chain ring
\end{theorem}

A {\it valuation ring} in a division ring $\Bbb{D}$ is a proper subring $R$ with the property that $x$ or $x^{-1}$ $\in R$ for each nonzero $x \in \Bbb{D}$. A ring is a valuation ring if and only if it is a left and right chain domain (i.e. a chain ring without zero divisors) (for a proof, see \cite{Lorim1}. Since any finite domain is a finite field, there do not exist finite valuation rings.

A ring $R$ is called  {\it E--ring} if and only if it possesses an ideal $I$ so that all ideals of $R$ are of the form $I^n$. In the infinite case $E$--rings can be characterized as $H$--rings with nilpotent radical and also as proper homomorphic images of discrete valuation rings (see \cite{Art}, \cite{As}, \cite{Lorim1} and \cite{Lorim2}). In the finite case the notions of $H$--ring and $E$--ring coincide.

The simplest and most investigated finite chain rings are the {\it Galois rings}, first defined by Krull \cite{Krull} as ``Grundringe" and later rediscovered by Janusz \cite{Janusz} and Raghavendran \cite{Rag}. A Galois ring is a commutative local PIR such that $J(R)=(p)$ with $p=1+1+ \ldots + 1$ ($p$ terms) for some prime $p$. 
These rings are very close to Galois fields. In the past ten years, finite chain rings and in particular Galois rings got a lot of attention in connection with coding theory (see e.g.~\cite{HoLa} and \cite{HoLa2}).

As for Galois fields one has

\begin{theorem} {\rm (Raghavendran \cite{Rag} and McDonald \cite{McD})}
For any prime $p \in \Bbb{N}$ and for any $n,r \in \Bbb{N}$ there exists a unique (up to isomorphism) Galois ring $R$ consisting of $q^n$ (with $q=p^r$) elements and with characteristic $p^n$
\end{theorem}

The unique Galois ring in the preceding theorem is denoted by GR($q^n,p^n$) (or sometimes also by GR($p^n,r$)). For $p=q$ we have  GR($p^n,p^n$) =  $\Bbb{Z}_{p^n}$, the ring of integers module $p^n$, and for $n=1$ we obtain the Galois field GF($q$)=GR($q,p$). 
All Galois rings can be constructed in the form $R=\Bbb{Z}_{p^n}[x]/(f(x))$ where $f(x)$ is a monic polynomial of degree $r$ which is irreducible modulo $p$ and hence GR($q^n,p^n)$, $q=p^r$, can be seen as Galois extensions of degree $r$ of its subring $\Bbb{Z}_{p^n}$. 

The properties of Galois rings are well known, e.g.~the structure of the group of units, the automorphism group, the possible subrings etc. Many results can be found in \cite{BiFl}.
The classification of all chain rings is still an open problem but partial results are known.
Galois rings occur in the construction of finite chain rings as can be seen from next theorem.

\begin{theorem}{\rm (Clark and Liang \cite{ClLi}, Wirt \cite{Wirt}, Neumaier \cite{Neu})}
Let $R$ be a finite chain ring with parameters $p, n, r$ and $k$ as in theorem 2.1. Then there exist integers $t$ and $s$ such that $R$ is isomorphic to $S[x,\sigma]/(g(x),p^{k-1}x^t)$ with $S$ = GR($q^k,p^k$) and $S[x,\sigma]$ the Ore--skew polynomial ring over $S$, i.e.~with usual addition and the multiplication $xa=\sigma(a)x$ for $\sigma \in Aut\; S$, and with $g(x)\in$ $S[x,\sigma]$ an Eisenstein polynomial of degree $s$, \mbox{$g(x)=x^s-p(a_0+a_1x + \ldots + a_{s-1}x^{s-1})$} with $a_0$ a unit in $S$ and $n=(k-1)s+t$ $\quad$($1 \leq t \leq s \leq n$)
\end{theorem}

The integer $s$ in the theorem above is called the ramification index of $R$. It is the smallest integer such that the ideal $(p)$ is equal to $J(R)^s$.

For a given set of parameters $p,n,r,k,s$ and $t$ one could ask for the number of non--isomorphic finite chain rings. In general this problem is still open, but partial results are known (see e.g.~\cite{Al}, \cite{Ark}, \cite{Ryb}, \cite{Nech} ). In some cases the parameters determine completely the ring:

\begin{theorem}{\rm (Clark and Liang \cite{ClLi} and Arkhipov \cite{Ark})}
Let $R$ be a finite chain ring with parameters $p,n,r,k,t$ and $s$ as in theorem 2.7. 

(a) If $k=1$ (hence $R$ has minimal characteristic $p$), then $R$ is uniquely determined (up to isomorphism) and $R \cong$ GF($q$)$[x,\sigma]/(x^n)$ (a truncated skew polynomial ring)

(b) If $k=n$ (hence $R$ has maximal characteristic $p^r$), then $R$ is uniquely determined (up to isomorphism) and $R \cong$ GR($q^n,p^n$) (always commutative)
\end{theorem}

Some more results are known for finite chain rings with characteristic $p^k$ with $1<k<n$ (see \cite{Hou}). In \cite{WYL} still another description of finite (commutative) chain rings is given as certain homomorphic images of the polynomial ring $\Bbb{Z}_{p^r}[x,y]$.

An important special subclass of chain rings is the one for which the Jacobson radical $J$ of $R$ has index of nilpotency $2$, so $J^2=0$. In this case $n=2$ and the two cases in theorem 2.8 are the only possible ones. Finite chain rings with $J \not=0$ and $J^2=0$ are called {\it uniform}.
The classification of finite uniform chain rings follows from theorem 2.8 but was also proved directly by Cronheim. 

\begin{theorem}{\rm (Cronheim \cite{Cron})}
Every finite uniform chain ring with $R/J\cong $ GF($q$) is either a ring of (twisted) dual numbers, or a truncated Witt ring of length 2, over the field GF($q$).
\end{theorem}

Rings of {\it (twisted) dual numbers} are the rings $\Bbb{D}(q,\sigma)$ = GF($q$)$[x,\sigma]/(x^2)$ (twisted for $\sigma \not= 1$) (corresponding to case ($a$) with $n=2$ in theorem 2.8).\newline
$\Bbb{D}(q,\sigma)$ can also be represented as the subring of matrices $\left( \begin{array}{ll} a & b \\ 0 & a^{\sigma} \end{array} \right)$ in the full matrix ring $M_2(q)$ of $2 \times 2$ matrices with elements in GF($q$).

$W_2(q)$, the {\it truncated Witt ring of length 2} over GF($q$) is defined on the set GF($q$) $\times$ GF($q$), $q=p^k$, as follows: 

addition: $(x_0,x_1)+(y_0,y_1)=(x_0+y_0,x_1+y_1+\dfrac{x_0^p+y_0^p-(x_0+y_0)^p}{p})$

and multiplication: $(x_0,x_1) \cdot (y_0,y_1)=(x_0y_0,x_0^py_1+x_1y_0^p)$

It can be proved that $W_2(q)$ is isomorphic to the Galois ring GR($q^2,p^2$) (this is case ($b$) with $n=2$ in theorem 2.8).

\section{Finite ring planes}

In this section we deal with geometries over finite rings and we restrict ourselves to the case of plane projective geometries. 
The projective line, higher dimensional projective geometries, affine and metric geometries, planar circle geometries (Benz--planes), chain geometries and polar geometries over finite rings will not be considered here.
 
A big part of finite geometry, called Galois geometry, is related to finite fields, see e.g.~the work of Hirschfeld \cite{Hir}. 
Since the pioneering work of B. Segre, the finite desarguesian (pappian) projective plane PG($2,q$) and its interesting point sets (arcs, ovals, blocking sets, unitals, $\ldots$) have been studied extensively. For planes over finite rings still a lot of work has to be done. As already mentioned in the introductory section, plane geometries over some important classes of rings have been defined in a suitable way, starting somewhere in 1940 by Barbilian (some isolated cases over particular rings were even known longer ago).

Before we look at planes over finite rings, we first recall the definition of a projective plane over an arbitrary ring (not necessarily finite). 

Let $R$ be an arbitrary ring (associative and with unit element). Denote the set of twosided invertible elements of $R$ by $R^{\star}$. Following \cite{Barb}, \cite{DeW} or \cite{Knu}, we can construct a plane projective geometry PG($2,R$) over $R$ as follows:

{\sl points} are the left unimodular triples $(x,y,z) \in R \times R \times R$ up to a right scalar in $R^{\star}$ (where $(x,y,z)$ left unimodular means that there exist $a,b,c \in R$ such that $ax+by+cz=1$ or equivalently $Rx+Ry+Rz=R$)

{\sl lines} are the right unimodular triples $[u,v,w] \in R \times R \times R$ up to a left scalar in $R^{\star}$ (where $[u,v,w]$ right unimodular means that there exist $a,b,c \in R$ such that $ua+vb+wc=1$ or equivalently $xR+yR+zR=R$)

{\sl incidence} I (between points and lines) is defined as follows: $(x,y,z)$ I $[u,v,w]$ if and only if $ux + v y + w  z=0$ 

{\sl neighborship} $\sim$ (between points and lines) is defined by:  $(x,y,z) \sim [u,v,w]$ if and only if  $u x + v  y + w z \in R\setminus R^{\star}$ \newline
It is clear that incidence always implies neighborship, so $p$ I $L$ implies $p \sim L$ for any point $p$ and any line $L$.

The incidence structure (with neighbor relation) obtained in this way is called {\it right projective plane over $R$}. In the same way one can define the {\it left projective plane over $R$} which is clearly isomorphic to the right projective plane over the opposite ring $R^{\circ}$. Therefore we will drop from now on the specification ``right" or ``left".

Although the denomination ``projective plane"  is used here, the projective plane over a ring (which is not a division ring) isn't a projective plane in the usual sense, as two distinct points may be incident with none or with more than one line and dually.

In addition to the neighbor relation for point--line pairs, also a neighbor relation between points (between lines respectively) can be considered in PG($2,R$):
points $(x,y,z)$ and $(x',y',z')$ are neighboring iff $\left[ \begin{array}{lll} x & y & z \\ x' & y' & z' \end{array} \right]$ cannot be extended to an invertible $3 \times 3$--matrix over $R$ (and similar for lines). 

Dealing with projective planes over rings it is  natural to assume that non--neighboring elements behave the same as distinct elements in an ordinary projective plane over a division ring.
To get that situation, it is necessary to restrict to the class of rings for which every one--sided unit is a two--sided unit, which was first observed by Barbilian \cite{Barb}.
Indeed, assume that $r$ is right--invertible (with right inverse $a$), but not left invertible. Consider the lines $[1,0,0]$ and $[r,0,0]$ (remark that $[r,0,0]$ is right unimodular since $r\cdot a+0\cdot b+ 0 \cdot c = 1$). These lines are distinct as otherwise there would exist a left scalar $l \in R^{\star}$ for which $[1,0,0]=l\cdot [r,0,0]$, so $1=l\cdot r$ in contradiction with the assumption that $r$ has not a left inverse. Now these two distinct lines are incident with the non--neighboring points $(0,1,0)$ and $(0,0,1)$.

The restriction to rings in which any left (or right) invertible element is a two--sided unit (or equivalently $a\cdot b=1$ implies $b\cdot a = 1$), the so--called Dedekind--finite rings,  only comes up when one deals with infinite rings. In the finite case (but also for many important classes of infinite rings) \mbox{invertible} elements are always two--sided invertible. In this context it is also interesting to mention that all reversible rings (i.e.~$a\cdot b = 0$ implies $b\cdot a =0$) are Dedekind--finite. 

Next we are interested in the connection between properties of the ring $R$ and the projective plane PG($2,R$). Most of the following results are reformulations for the finite case of theorems that can be found in Veldkamp \cite{Veld1} and \cite{Veld2} .

In the projective plane over a Dedekind--finite ring, there is a unique line incident with two given non--neighboring points (and dually). One might ask whether the neighbor relation is completely determined by the incidence relation in the sense that two points are non--neighboring if and only if there is unique line incident with them (and dually for lines). This is not always the case, but it does if every non--invertible element in $R$ is a right and left zero--divisor.

A (left or right) artinian ring is a ring in which any non--empty set of (left or right) ideals that is partially ordered by inclusion, has a minimal element. In a (left or right) artinian ring any non--invertible element is a (left or right) zero divisor. As finite rings are always left and right artinian we get that the neighbor relation is completely determined by the incidence relation in planes PG($2,R$) over finite rings. In \cite{Corb0} a proof is given of the property that in a finite ring any left zero divisor is also a right zero divisor.

\begin{theorem}{\rm (Veldkamp \cite{Veld1})}
Let $R$ be a finite ring and PG($2,R$) the projective plane over $R$. Then two distinct points are neighboring if and only if they are incident with either no or at least two lines. Two distinct  lines are neighboring if and only if they are incident with either no or at least two points.
\end{theorem} 

A projective ring plane is called {\it linearly connected} (Veldkamp \cite{Veld1}), {\it neighbor cohesive} (Drake and Jungnickel \cite{DrJ}) or {\it punctally cohesive} (Baker et al. \cite{LO3}) if any two distinct points are incident with at least one line.
For planes over rings of stable rank 2 it is proved in Veldkamp that two points are incident with at least one line if and only if $R$ has the following property: for any two $r_1, r_2 \in R$ there exists $a \in R$ such that $Rr_1+Rr_2=R(r_2+ar_1)$. This is fulfilled for $R$ a left B\'ezout ring, i.e. a ring for which any finitely generated left ideal is a principal ideal. Dually two lines are incident with at least one point iff $R$ is a right B\'ezout ring.
For finite rings the B\'ezout conditions amount to the condition that $R$ is a principal ideal ring (recall that for a finite ring the notions left principal, right principal and principal coincide).
So we can reformulate the theorem for finite rings as follows:

\begin{theorem}{\rm (Veldkamp \cite{Veld1})}
Let $R$ be a finite ring and PG($2,R$) the projective plane over $R$. Then any two points are incident with at least one line (the plane is linearly connected) and dually if and only if $R$ is a principal ideal ring.
\end{theorem}

The possibility of more than one line incident with two neighboring points (and dually) corresponds to the presence of zero divisors in the ring. So any two distinct points are incident with exactly one line (and dually) if and only if $R$ is a B\'ezout domain. In the finite case this becomes: if $R$ is a principal ideal domain, hence if $R$ is a finite field. Hence:

\begin{theorem}{\rm (Veldkamp \cite{Veld1})}
Let $R$ be a finite ring and PG($2,R$) the projective plane over $R$. Then any two distinct points are incident with exactly one line and dually if and only if $R$ is a finite field (i.e.~PG($2,R$) is a pappian projective plane).
\end{theorem}

Next we look at the special case of (finite) {\it local} rings. For such rings the definition of the projective plane PG($2,R$) and his neighbor relations, can be adapted (in an equivalent way) a little. E.g.~two points $(x,y,z)$ and $(x',y',z')$ are neighbors if and only if $(x',y',z')$ $-$ $(x,y,z)\lambda \in J \times J \times J$ for
some $\lambda \in R \setminus J$ with $J$ the maximal ideal of $R$ and similarly for lines.

\begin{theorem}{\rm (Veldkamp \cite{Veld1})}
Let $R$ be a (finite) ring and PG($2,R$) the projective plane over $R$. Then the neighbor relation $\approx$ between points (between lines resp.) is transitive if and only if $R$ is a (finite) local ring.
\end{theorem}

For local rings $R$ there is a canonical epimorphism $\varphi$ from $R$ onto the division ring $\Bbb{K}=R/J$. This epimorphism induces an epimorphism $\pi$ of the projective plane PG($2,R$)  onto the (ordinary) projective plane PG($2,\Bbb{K}$) by putting $\pi(x,y,z)=(\phi(x),\phi(y),\phi(z))$ and $\pi[u,v,w]=[\phi(u),\phi(v),\phi(w)]$ and the neighbor relation can be expressed by means of $\pi$ : $p \sim L$ if and only if $\pi(p)$ I $\pi(L)$ and similarly $p \sim q$ iff $\pi(p) = \pi(q)$ and $L \sim M$ iff $\pi(L)=\pi(M)$. 

The projective plane over a local ring, also known as a (desarguesian) projective Klingenberg plane, therefore is strongly connected with an ordinary desarguesian projective plane. One could say that the points (and lines) of an ordinary projective plane are blown up to clusters of neighboring points (lines) to produce a projective Klingenberg plane. In the finite case the epimorphic image of PG($2,R$) is the plane PG($2,q$) over the Galois field GF($q$).

Combining theorems 3.1, 3.2 and 3.4 yields the following:

\begin{theorem}{\rm (Veldkamp \cite{Veld1})}
Let $R$ be a finite ring and PG($2,R$) the projective plane over $R$. Then the neighbor relation $\approx$ between points (between lines resp.) is transitive and two neighboring points are incident with at least two lines and dually, if and only if $R$ is a finite local principal ideal ring.
\end{theorem}

From section 2 we know that finite local principal ideal rings are synonym for finite chain rings or finite $H$--rings.
Recall that projective planes over $H$--rings are called (desarguesian) projective Hjelmslev planes. 

We now summarize the possibilities for projective planes over a finite ring.

\begin{corollary}
Let $R$ be a finite ring and PG($2,R$) the projective plane over $R$. Then only four cases are possible:\newline
(a) $R$ has no zero divisors and hence is a field and PG($2,R$) is an ordinary pappian projective plane (two distinct points are incident with exactly one line and dually)\newline
(b) $R$ is a local principal ideal ring (hence a chain ring = $H$--ring) and PG($2,R$) is a desarguesian projective Hjelmslev plane \newline
(c) $R$ is local but not a principal ideal ring and PG($2,R$) is a desarguesian projective Klingenberg plane (but not a Hjelmslev plane)\newline
(d) $R$ is semilocal (but not local) and PG($2,R$) has non--transitive neighbor relations
\end{corollary}

The fourth class (d) is the wildest as it contains all finite rings which are not local (but necessarily semilocal due to the finiteness, i.e.~with a finite number of maximal ideals). Important examples of rings belonging to this class are the full matrix rings $M_n(q)$ over GF($q$). Projective planes over full matrix rings were first mentioned by Ree \cite{Ree} and further studied by J.A.~Thas \cite{Thas}  who also gave an interpretation of PG($2,M_n(q)$) in terms of the projective space PG($3n-1,q$). 
Other examples are the rings $\Bbb{Z}_m$ with $m\not=p^r$ (see \cite{EuGa}). Of special interest are also the rings of double numbers $\Bbb{B}(q)=$GF($q$) + GF($q$)$\,t$ with $t^2=t$. They possess exactly two maximal ideals. In \cite{SaPl} projective planes over $\Bbb{B}(q)$ are studied.

Examples of finite local rings that are not chain rings (c) are provided by the rings GF($q$)[$x,y$]/$\langle x^n,xy,y^n \rangle$ ($n>1$). The corresponding planes are finite desarguesian Klingenberg planes that are not Hjelmslev planes.

Class (b) contains many interesting examples, including the Galois rings GR($q^n,p^n$) ($q=p^r$) and the rings $\Bbb{A}(p^r,n)=$ GF($p^r$)$[x]/x^n$ (called quasi--Galois rings in \cite{BiFl}). The rings $\Bbb{A}(p^r,n)$ can also be interpreted as matrix rings, consisting of all matrices $(a_{ij})$ with elements belonging to GF($q$) and $a_{i,j+i-1}=a_{1j}$ and $a_{ij}=0$ for $i>j$.
For $n=2$ the ring of dual numbers $\Bbb{D}(q)$ over GF($q$) is included. Projective planes over dual numbers were considered yet a century ago by Corrado Segre \cite{Segre}. 

Class (a) finally consists of all the Galois fields GF($q$) with the associated projective planes PG($2,q$).

For finite (not necessarily desarguesian) projective Klingenberg and Hjelmslev planes a unique set of para\-meters (the order) can be given (see \cite{Klein} and \cite{DrL}): for any flag $(p,L)$ there
are exactly $t$ points on $L$ neighboring with $p$ and exactly $s$
points on $L$ not neighboring with $p$. Moreover: the number of points = the number of lines = $s^2+st+t^2$, any line is incident
with $s+t$ points, any point is incident with $s+t$ lines, any point
has $t^2$ neighbors, any line has $t^2$ neighbors, $t|s$ and
$r=\frac{s}{t}$ is the order of the  projective plane that is the canonical epimorphic image of the Klingenberg plane, and $s\leq t^2$ or $t=1$.

For a finite desarguesian projective Klingenberg plane this yields  

\begin{theorem} {\rm (Drake and Jungnickel \cite{DrJ} )}
Let $R$ be a finite local ring. Then the projective Klingenberg plane (Hjelmslev plane in some cases) PG($2,R$) has parameters $s=|R|=q^n$ and
$t=|J|=q^{n-1}$ with $q=p^r$ a prime power.

\end{theorem} 

To conclude this section we consider rings of order 4. This is the smallest order for which there exist rings that are not division rings. There are four non--isomorphic rings (with unit) of order 4. The first is the Galois ring GF($4$) that gives rise to the projective plane PG($2,4$). The second is the chain ring $\Bbb{Z}_4 \cong GR(4,4) \cong W_2(2)$ with characteristic 4, the coordinate ring of a projective Hjelmslev plane. The third is the chain ring $\Bbb{D}(2) \cong \Bbb{A}(2,2)$ of characteristic 2 and with GF($2$) as a subfield (dual numbers over GF($2$)) that gives rise to another projective Hjelmslev plane (it is proved in \cite{KepVM} that the plane over $\Bbb{D}(2)$ is embeddable in PG($5,q$) while the plane over $\Bbb{Z}_4$ is not). Finally the fourth is the non--local ring $\Bbb{B}(2) \cong GF(2)[t]/(t^2-t)$ of characteristic 2 (double numbers over GF($2$)), associated to a Veldkamp plane with non--transitive neighbor relation.

\section{Finite ring--like structures and ring--like planes}

Besides finite rings, even more general finite ring--like algebraic structures deserve a closer look in relationship to geometry. Till now very few research has been done in this area (except for finite field--like algebraic structures).\newline
In the literature several generalizations of rings can be found. Among the most important are: non--associative rings, nearrings and semirings. For all these structures there are finite examples and a generalization of the concept ``local ring" exists, which opens perspectives for Klingenberg--like geometries associated to these generalized rings. 

\subsection{Semirings}
A {\sl semiring} is a structure $(S,+,\cdot)$ with $(S,+)$ a commutative semigroup with identity element $0$, $(S,\cdot)$ a (not necessarily commutative) semigroup with identity element $1$ ($\not=0$), in which the left and right distributivity of multiplication over addition hold and in which $0$ is absorbing for multiplication $a\cdot 0=0 \cdot a=0$. Hence semirings differ from rings by the fact that elements do not always have an inverse for the addition (the additive group of a ring is replaced by a semigroup).\newline
Semirings were first introduced by Vandiver \cite{Vand} in 1935 and in the past years there was an enormous amount of publications on the subject (see e.g. the work of Glazec \cite{Gla} for a survey), mainly in relation to computer science and automata theory but they also are interesting algebraic objects on their own, see \cite{HW} and \cite{Gol}.

Examples of finite semirings are $B(n,i)$ on the set $\{0,1,\ldots,n-1\}$ ($n \geq 2$ and $0 \leq i \leq n-1$) with addition $\oplus$ defined by $a\oplus b =a+b$ if $0 \leq a+b < n$ and $a \oplus b=c$ if $a+b \geq n$ with $c$ the unique number such that $c \equiv a+b$ (mod $n-i$) and $0\leq c \leq n-1$. Multiplication $\odot$ is defined in a similar way. In particular $B(n,0)$ is the ring $\Bbb{Z}_n$ of integers modulo $n$ and $B(2,1)$ is known as the boolean semiring $\Bbb{B}$. For other values of $n$ and $i$ one obtains semirings that are not rings. 

In semirings zero divisors and zero sums are of interest ($a$ is a zero sum of $S$ if there exists an element $b\not=0$ in $S$ such that $a + b = 0 $).
In \cite{HL} it is proved that if $S$ is a finite commutative semi\-ring then either \mbox{every} zero sum is a zero divisor or $S$ is a ring. As a corollary one has that a finite commutative semiring without zero divisors (a semidomain) either is zero sum free ($a+b=0$ always implies $a=b=0$) or is a finite domain (and hence a field).

A semiring is called a semifield if $(S^{\ast},\cdot)$, with $S^{\ast}=S \setminus \{0\}$, is a group
(in a semifield the multiplication need not to be commutative, so the term semi division ring would be better). We have to warn for confusion between the semifields considered here and the semifields defined in the context of non--desarguesian projective planes by e.g. Albert, Dickson, Knuth and others and which are also known under the name division algebras or distributive quasifields. Those are generalizations of division rings by dropping the need for associativity of the multiplication, so $(S,+,\cdot)$ is in that context a semifield if $(S,+)$ is a group (which turns out to be commutative) and $(S^{\ast},\cdot)$ is a loop and the two distributivity laws of multiplication over addition hold (for a survey on those semifields, see e.g.~\cite{Lav})

To my knowledge no research has been done yet on geometries over (finite) semirings that are not rings. So there may be opportunities. In particular a generalization to Klingenberg--like planes (over local semirings) seems possible though not trivial as the concept of ideals in semirings is much more complicated. As a starting point for ideals in semirings and the concept of local semiring see e.g. \cite{At} and \cite{Gu}.

\subsection{Nearrings}
A {\sl left nearring} is a structure $(N,+,\cdot)$ with $(N,+)$ a (not necessarily) commutative group with identity element $0$, $(N,\cdot)$ a (not necessarily commutative) semigroup with identity element $1$ ($\not=0$), in which the left distributivity of multiplication over addition holds: $a\cdot (b+c)=a\cdot b + a \cdot c$. Similar for {\sl right nearring}. Nearrings differ from rings by the fact that addition isn't necessarily commutative and there is only distributivity on one side. Nearrings which are distributive on both sides are rings (the commutativity of the addition then follows automa\-tically). Most of the material on nearrings can be found in the work of Pilz \cite{Pilz1} and \cite{Pilz2}.

A (left or right) nearring is called a (left or right) nearfield if $(S^{\ast},\cdot)$ is a group, with $S^{\ast}=S \setminus \{0\}$. So nearfields are division rings with only distributivity on one side. Nearfields were first discovered by Dickson in 1905 and are useful in constructing examples of non--desarguesian projective planes. All finite nearfields are classified by Zassenhaus. They are either Dickson--nearfields or they belong to one of seven exceptional classes. 

Little research has been done on geometry over nearrings that are not near\-fields (except for the class of planar nearrings which give rise to balanced incomplete block designs, but these nearrings don't possess a multiplicative identity element and therefore are less usable in the context of projective geometry). So another suggestion for research could be a treatment of plane projective geometries over (finite) nearrings that are not rings. The special case of Klingenberg planes (over local nearrings) and Hjelmslev planes (over $H$--nearrings) was already initiated in two general papers by E. Kolb, see \cite{Kolb1} and \cite{Kolb2} and by T\"orner in \cite{bibitem65}.
Local nearrings were introduced by Maxson in \cite{Maxson2} and partially classified in \cite{Maxson3}.  Among the results is the fact that the additive group of a finite local nearring always is a $p$--group and the existence of a natural epimorphism from a local nearring onto a local nearfield. 
Other results on finite nearrings can be found in \cite{AHS}, \cite{Clay1},  \cite{Clay2}, \cite{Clay3}, \cite{Hub}, \cite{Ligh}, \cite{LighMa}, \cite{Maxson1} and \cite{Wendt}.

\subsection{Alternative rings}
Finally, we mention some results on non--associative rings.
A {\it non--associative ring} is a structure $(A,+,\cdot)$ which satisfies all axioms for an (associative) ring with multiplicative identity element, except for the associativity of the multiplication.
An {\it alternative ring} is a non--associative ring such that $a\cdot (a \cdot b)=a^2\cdot b$ and $(a\cdot b)\cdot b = a \cdot b^2$ for all $a,b \in A$. Alternativity is a weaker condition than associativity. 
If any element in an alternative ring $A$ is a unit, then $A$ is called an {\it alternative division ring}. Alternative division rings are used to construct a class of non--desarguesian projective planes, called Moufang planes. By the theorem of Artin--Zorn every finite alternative division ring is a field, hence finite Moufang projective planes are desarguesian (and pappian).

Generalizations to alternative rings that are not division rings are due to Baker, Lorimer and Lane. In \cite{Baker2} Moufang projective Klingenberg planes are defined as projective Klingenberg planes that are $(p,L)$--transitive for all flags and it is proved that they can be coordinatized by a local alternative ring. In \cite{Baker1}  several characterizations of local alternative rings are given and an analogue for non--associative chain rings and $H$--rings is defined properly. In the finite case it is proved that the concepts of alternative $H$--ring, left (or right) alternative chain ring and local alternative principal ideal ring are equivalent. Moreover the theorem of Artin--Zorn is expanded : any finite alternative chain ring (or $H$--ring) is associative \cite{LO3}.\newline  
Leaving out the condition of being local, leads to more general alternative rings. It is hard to define projective ring planes over such rings in a suitable way. Faulkner has done it for alternative stable rank 2 rings in \cite{Fauk1} (generalizing Veldkamp's results for associative stable rank 2 rings) and for alternative rings in which any one--sided unit is two--sided in \cite{Fauk2} (generalizing the planes of Barbilian).

\end{document}